\documentclass{amsart}

\usepackage[dvips]{graphicx}
\newtheorem{theorem}{Theorem}[section]
\newtheorem{proposition}[theorem]{Proposition}
\newtheorem{lemma}[theorem]{Lemma}
\newtheorem{claim}[theorem]{Claim}

\newtheorem{corollary}[theorem]{Corollary}

\theoremstyle{definition}

\newtheorem{example}[theorem]{Example}

\theoremstyle{remark}
\newtheorem{remark}[theorem]{Remark}

\numberwithin{equation}{section}

\begin{document}

\title[Crosscap numbers of torus knots]
{Crosscap numbers of torus knots}
\author{Masakazu Teragaito}
\address{Department of Mathematics and Mathematics Education,
Faculty of Education, 
Hiroshima University, 1-1-1 Kagamiyama, Higashi-hiroshima 739-8524, Japan}
\email{teragai@hiroshima-u.ac.jp}



\begin{abstract}
The crosscap number of a knot in the $3$-sphere is the minimal genus of non-orientable
surface bounded by the knot.
We determine the crosscap numbers of torus knots.
\end{abstract}
\maketitle

\section{Introduction}\label{intro}

For a knot $K$ in the $3$-sphere $S^3$, the \textit{crosscap number\/} of $K$, denoted by $c(K)$,
is defined to be the minimal first betti number of compact, connected, non-orientable
surfaces $F$ with $\partial F=K$ \cite{C}.
(For the trivial knot, it is defined to be zero instead of one.)
Since the crosscap number is an analogy of the genus of a knot, it is also called 
the non-orientable genus \cite{MY2}.
Clark \cite{C} showed that a knot has crosscap number one if and only if
it is a $2$-cable of some knot.
In \cite{T1}, we showed that genus one, crosscap number two knots are doubled knots.
Since a genus one knot has crosscap number at most three (\cite{C,MY}), this determines the crosscap numbers of genus one knots.
However it is hard to determine the crosscap number of a knot, in general.  See also \cite{B,IOT,MY,U}.

Let $F$ be a surface bounded by a knot $K$.
Then it can be assumed that $F$ meets the regular neighborhood $N(K)$ of $K$
in a collar neighborhood of $\partial F$ in $F$, which is an annulus.
Then $F\cap \partial N(K)$ is an essential simple closed curve on $\partial N(K)$.
Its (unoriented) isotopy class is referred to as the \textit{boundary slope\/} of $F$, which
is parameterized by integers in the usual way (see \cite{R}).
The boundary slope of $F$ is uniquely determined by $F$.
Clearly, if $F$ is orientable, then the boundary slope is zero.
But this is not the case when $F$ is non-orientable.
Then the boundary slope must be an even integer by homological reason
(see \cite{BW}).
If $F$ is a spanning surface of $K$, that is $\partial F=K$, then 
a new non-orientable spanning surface $F'$ of $K$ is obtained by
adding a small half-twisted band to $F$ locally.
The boundary slope of $F'$ is $(\text{that of}\ F)\pm 2$.
Thus any even integer can be the boundary slope of some non-orientable spanning surface of $K$.
The first betti number of $F$ is denoted by $\beta_1(F)$.
In particular, $\beta_1(F)$ is equal to the genus of $F$
when $F$ is non-orientable. 

In this paper, we determine the crosscap number of a torus knot $K$ and the boundary slope
of a non-orientable surface bounded by $K$, which realizes its crosscap number $c(K)$, simultaneously.
In fact, we show that the boundary slope of a non-orientable spanning surface of $K$ realizing
$c(K)$ is uniquely determined by $K$.
This does not hold in general.  For example, the figure-eight knot has crosscap number two, and
it bounds two once-punctured Klein bottles with boundary slopes $4$ and $-4$,
which are two checkerboard surfaces of a standard diagram.
Another remarkable example is the $(-2,3,7)$-pretzel knot, which bounds two once-punctured
Klein bottles with boundary slopes $16$ and $20$. See \cite{IOT}. 
For torus knots, the crosscap numbers are given by using a positive integer $N(x,y)$ introduced by Bredon and Wood \cite{BW},
which is the minimal genus of closed, connected, non-orientable surfaces contained in the lens space $L(x,y)$.
It is easy to calculate $N(x,y)$ by using continued fractions or a recursive formula (see Section \ref{ex}).

Let $T$ be a standard torus in $S^3$.
It decomposes $S^3$ into two solid tori $V$ and $W$.
Let $f:S^1\times D^2\to V$ be a homeomorphism.
This determines the longitude-meridian system of $V$ by setting
$\lambda=f(S^1\times *), *\in \partial D^2$, and $\mu=f(*\times \partial D^2), *\in S^1$, which
gives a basis of $H_1(T)$.
We assume that $\lambda$ is a preferred longitude, that is, $\lambda$ bounds a disk in $W$.
Let $T(p,q)$ be the torus knot of type $(p,q)$ lying on $T$,
which represents $p[\lambda]+q[\mu]$ in $H_1(T)$.
Note that $p$ and $q$ are coprime.
Since all $T(p,q),T(-p,q),T(p,-q), T(-p,-q), T(q,p)$ are equivalent (there is a homeomorphism of $S^3$
sending one to the other), 
they have the same crosscap number. Therefore we always assume $p,q>0$.  
The torus knot $T(p,q)$ is said to be \textit{odd\/} (resp.~\textit{even\/})
if $pq$ is odd (resp.~even).
Furthermore, if $T(p,q)$ is odd (resp.~even), then we always assume that $p>q$
(resp.~$p$ is even).

\begin{theorem}\label{main}
Let $K$ be the non-trivial torus knot of type $(p,q)$, where $p,q>0$, and
let $F$ be a non-orientable spanning surface of $K$ with $\beta_1(F)=c(K)$.
\begin{itemize}
\item[(1)] If $K$ is even, then $c(K)=N(p,q)$ and
the boundary slope of $F$ is $pq$.
\item[(2)] If $K$ is odd, then
$c(K)=N(pq-1,p^2)$ \textup{(}resp.~$N(pq+1,p^2)$\textup{)} and
the boundary slope of $F$ is $pq-1$ \textup{(}resp.~$pq+1$\textup{)}
if $xq\equiv -1\pmod p$ has an even \textup{(}resp.~odd\textup{)} solution $x$\ satisfying $0<x<p$.
\end{itemize}
\end{theorem}

Remark that the equation $xq\equiv -1\pmod p$ in (2) has the unique solution $x$ satisfying $0<x<p$.
See Section \ref{odd}.

In general, the crosscap number is not additive under the connected sum operation \cite{MY}.
But we have:

\begin{theorem}\label{main2}
If $K_1,K_2,\dots,K_n$ are torus knots, then
$$c(K_1\sharp K_2\sharp \dots \sharp K_n)=c(K_1)+c(K_2)+\dots +c(K_n).$$
\end{theorem}


\section{Preliminaries}\label{pre}

Let $K=T(p,q)$ be the non-trivial torus knot of type $(p,q)$, and let $E(K)$ be its exterior.
Let $F$ be a non-orientable surface bounded by $K$ realizing its crosscap number.
We may assume that $F\cap N(K)$ is an annulus, and $F\cap E(K)$ is also denoted by $F$.
As noted in Section \ref{intro}, $\partial F\ (\subset \partial E(K))$ determines
the boundary slope $r$, which is an even integer.

\begin{lemma}\label{incompressible}
$F$ is incompressible in $E(K)$.
\end{lemma}

\begin{proof}
Assume that $F$ is compressible in $E(K)$, and let $D$ be a compressing disk for $F$.
Note that $\partial D$ is orientation-preserving in $F$.

If $\partial D$ is separating in $F$, then
compression along $D$ gives two surfaces $F_1$ and $F_2$, where
$\partial F_1=\partial F$ and $F_2$ is closed.
Then $F_2$ is orientable, and hence $F_1$ is non-orientable.
It is easy to see that $\beta_1(F_1)\le \beta_1(F)-2$.
This contradicts the minimality of $\beta_1(F)$.
Therefore $\partial D$ is non-separating in $F$.

Let $F'$ be the resulting surface obtained by compressing $F$ along $D$.
Since $\beta_1(F')=\beta_1(F)-2$, $F'$ must be orientable by the minimality of $F$.
But if we add a small half-twisted band to $F'$ (after extending $F'$ to $K$ radially in $N(K)$),
then we obtain a non-orientable surface $F''$ bounded $K$ with $\beta_1(F'')=\beta_1(F')+1<\beta_1(F)$.
This contradicts the minimality of $F$.
\end{proof}

Recall that $K$ lies on the standard torus $T$.
Let $A=T\cap E(K)$, and let $a_1$ and $a_2$ be the components of $\partial A$.
Then it is well known that $A$ is an essential (incompressible and boundary-incompressible)
annulus in $E(K)$ and that each component $a_i$ of $\partial A$ has slope $pq$ on $\partial E(K)$.
By an isotopy of $F$, we may assume that $F$ and  $A$ intersect transversely, and therefore $F\cap A$ consists
of arcs and loops.
Since $F$ and $A$ are incompressible, we can remove any loop component of $F\cap A$ bounding a disk on
either $F$ or $A$.
Furthermore, we may assume that $\partial F$ meets each $a_i$ in the same direction
(after giving them orientations).
Then $F\cap A$ contains exactly $\Delta(r,pq)=|r-pq|$ arcs, since $\partial F$ meets $a_i$ in
$\Delta(r,pq)$ points, where $\Delta(r,pq)$ denotes the minimal geometric intersection number between
two slopes $r$ and $pq$ on $\partial N(K)$.

\begin{lemma}\label{span}
If $F\cap A$ contains an arc component $\alpha$, then $\alpha$ connects distinct boundary components of $A$.
\end{lemma}

\begin{proof}
Suppose that the ends of $\alpha$ lie in $a_1$, say. 
Without loss of generality, we can assume that $\alpha$ is outermost on $A$.
That is, $\alpha$ cuts off a disk $D$ from $A$ such that $\mathrm{Int}D\cap F=\emptyset$.
Let $\partial D=\alpha\cup \beta$, where $\beta$ is a subarc of $a_1$.
See Figure \ref{compress}.

\begin{figure}[tb]
\includegraphics*[scale=0.7]{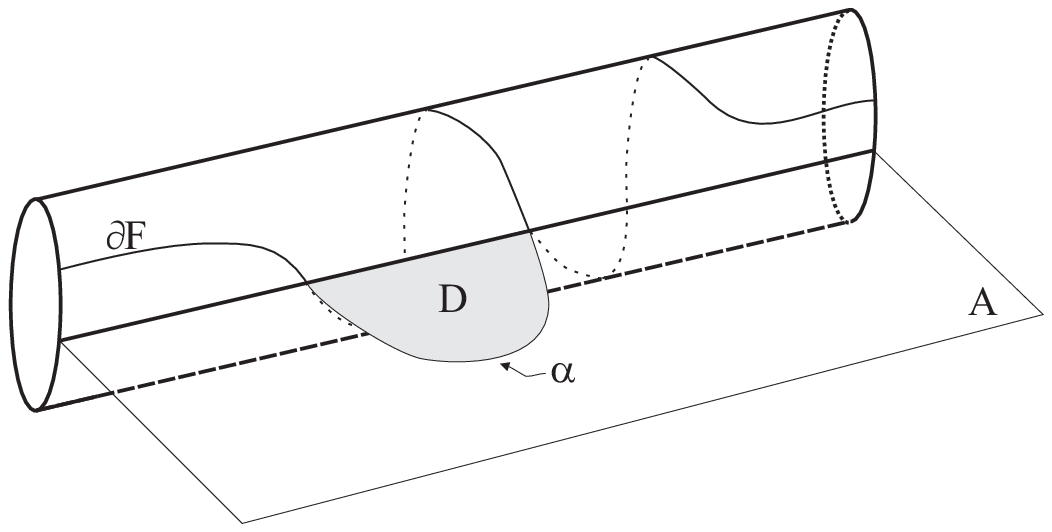}
\caption{}\label{compress}
\end{figure}

By doing boundary-compression along $D$, we obtain a connected surface $F'$.
Note that $\beta_1(F')=\beta_1(F)-1$ by an Euler characteristic calculation.
Since $\partial F'$ still meets a meridian of $K$ in one point, $F'$ can be extended to a spanning surface of $K$.
Thus $F'$ is orientable, because of the minimality of $F$.
Then $F$ is obtained from a Seifert surface $F'$ of $K$ by adding a half-twisted band locally.
Since $K$ has genus $(p-1)(q-1)/2$ (\cite{R}), it implies that $c(K)=\beta_1(F)=(p-1)(q-1)+1$. 
By \cite[Proposition 1.3]{MY}, $c(K)\le \min\{(p-1)q/2,(q-1)p/2\}$.
(A standard diagram of $K$ has $(p-1)q$ or $(q-1)p$ crossings.)
This gives a contradiction easily.
\end{proof}

Thus, if $|pq-r|\ne 0$, then $F\cap A$ contains no loop component.
If $|pq-r|=0$, then $F\cap A$ may contain loop components which are essential in $A$.
Since $r$ is even, $|pq-r|=0$ happens only when $K$ is even.
Also, if $K$ is odd, then $F\cap A$ contains at least one arc.

\begin{figure}[tb]
\includegraphics*[scale=0.7]{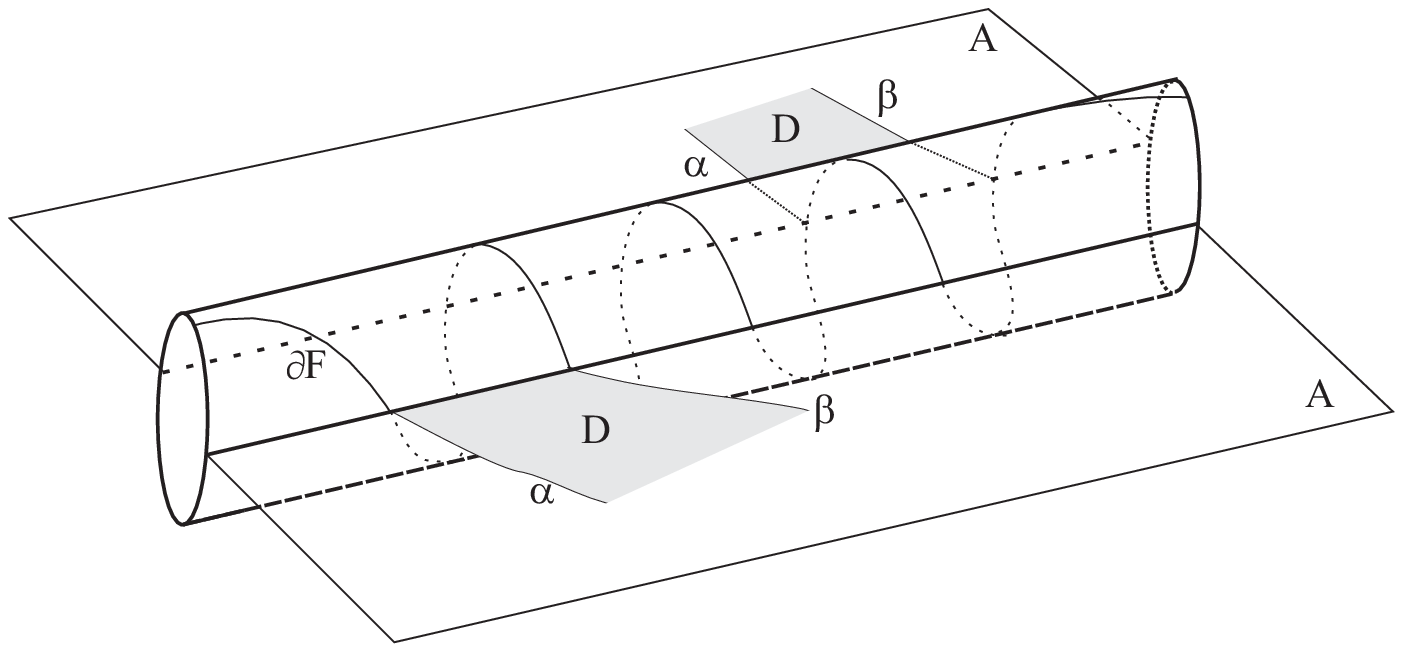}
\caption{}\label{split}
\end{figure}

We introduce an operation called a \textit{disk splitting}.
This operation will be used to determine the boundary slope of $F$ in the following sections. 
If $F\cap A$ contains at least two arcs, then there are
two arcs $\alpha$ and $\beta$ of $F\cap A$, which
cut a rectangle $D$ from $A$ such that $\mathrm{Int}{D}\cap F=\emptyset$. See Figure \ref{split}.
Let $F'$ be the surface obtained by splitting $F$ along $D$.
That is, take a product neighborhood $D\times [0,1]$ of $D$ in $E(K)$, and
let $F'=(F-D\times [0,1])\cup D\times\{0,1\}$.
Note that this disk splitting does not change Euler characteristic.
(If $F'$ is disconnected, then $\chi(F')$ is the sum of Euler characteristic of its components.)

\section{The case where $K$ is even}\label{even}

In this section, we prove Theorem \ref{main}(1).

Suppose that $K$ is even.
Recall that $p$ is assumed to be even, and that
$r$ is the boundary slope of $F$.

\begin{proposition}\label{reduction-even} 
$r=pq$.
\end{proposition}

\begin{proof}
Suppose not. 
Let $\Delta=\Delta(r,pq)=|r-pq|$.
Then $\Delta$ is even and $\Delta\ge 2$.

\begin{claim}\label{p4}
$p\ge 4$.
\end{claim}

\begin{proof}[Proof of Claim \ref{p4}]
If $p=2$, then $c(K)=1$ and $F$ is a M\"obius band.
Consider $r$-Dehn surgery $K(r)$.
That is, $K(r)$ is the union of $E(K)$ and a solid torus $J$
glued to $E(K)$ along their boundaries in such a way that $r$ bounds a meridian disk in $J$.
Thus $K(r)$ contains a projective plane, which is obtained by capping the M\"obius band
off by a meridian disk of $J$.
Then $K(r)$ is $P^3$ or reducible.
By \cite{M}, the former is impossible, and so $K(r)$ is reducible.
Also $pq$ is the only slope yielding a reducible manifold.
Hence $p\ge 4$.
\end{proof}

Now, $F\cap A$ consists of $\Delta$ arcs.
There are mutually disjoint $\Delta/2$ rectangles on $A$ cut by these arcs.
We perform $\Delta/2$ disk splittings to obtain two surfaces $F_1$ and $F_2$.
Note that $\chi(F)=\chi(F_1)+\chi(F_2)$.
One surface is connected and has a single boundary component, whose slope is $pq$, and the other has $\Delta/2$ 
boundary components, each of which is inessential on $\partial E(K)$.
In particular, the latter does not contain a non-orientable component.
See Figure \ref{splitting} (where $\Delta=6$).
We can assume that $F_1$ has a single boundary component.
Since the slope $pq$ of $\partial F_1$ is not zero, $F_1$ must be non-orientable.
If we show $\chi(F_2)<0$, then $\chi(F)<\chi(F_1)$, which contradicts the minimality of $F$.

\begin{figure}[tb]
\includegraphics*[scale=0.46]{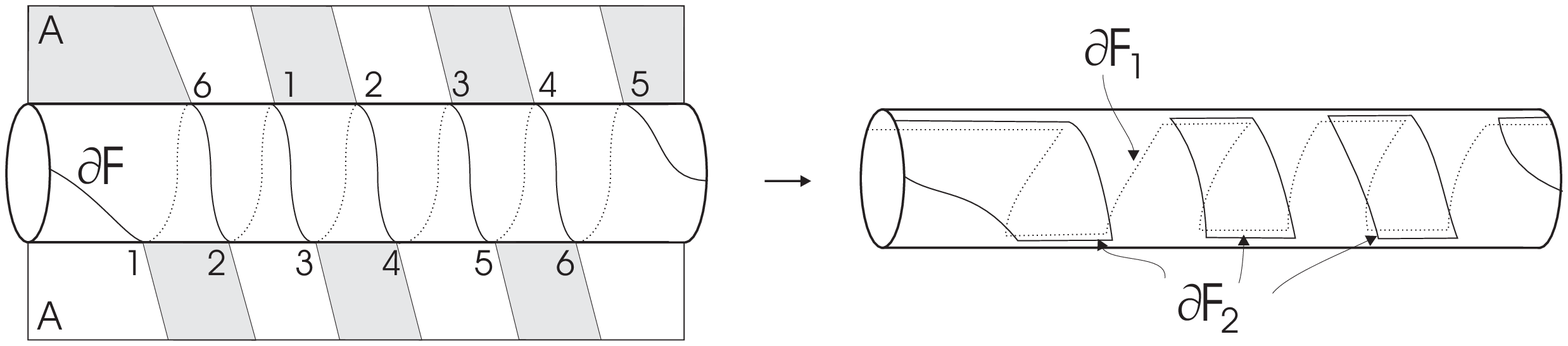}
\caption{}\label{splitting}
\end{figure}

Number the points of $a_1\cap \partial F$ $1,2,\dots,\Delta$ along $a_1$.
This induces the labeling of the points of $a_2\cap \partial F$ through the arcs of $F\cap A$.
See Figure \ref{splitting}.
The annulus $A$ splits $E(K)$ into two solid tori $U_1$ and $U_2$.
We assume $U_i$ contains $F_i$ for $i=1,2$.
Let $A_i=\partial E(K)\cap U_i$.
Then for some permutation $\sigma$ of $\{1,2,\dots,\Delta\}$,
the point with label $i$ in $a_1$ is connected to one with label $\sigma(i)$ in $a_2$
by the arcs of $A_2\cap \partial F$.
(In Figure \ref{splitting}, $\sigma=(153)(264)$.)
Note that $\sigma$ has an even number of orbits, all of which have the same length $\ell(\sigma)$, say, by the definition.
(Thus $\sigma$ has $\Delta/\ell(\sigma)$ orbits.)

Since $F_2$ has $\Delta/2$ boundary components, and $F_2$ has no closed components,
$|F_2|\le \Delta/2$.  ($|\dots|$ denotes the number of components.) 
In fact, the boundary components of $F_2$ belong to the same component of $F_2$ 
when they correspond to the same orbit of $\sigma$. 
Thus $|F_2|\le \Delta/2\ell(\sigma)$.

\begin{claim}\label{F2}
$F_2$ contains no disk components.
\end{claim}

\begin{proof}[Proof of Claim \ref{F2}]
Suppose not.
We can regard $F_2$ as the union of $F\cap U_2$ and the rectangles on $A$ used in the disk splittings.
Let $D$ be a disk component of $F_2$.
Then we have $\ell(\sigma)=1$.
We can take a loop $f$ on $\partial U_2$, which is the union of an arc of $F\cap A$ and an arc of $\partial F\cap A_2$
and which meets $\partial D$.
Then this loop $f$ bounds a disk in $D$.
Since $f$ meets the core of $A_2$ in one point, $f$ must be a meridian of $U_2$.
But this implies that $K$ is trivial, a contradiction.
\end{proof}

Thus $\chi(F_2)\le 0$.
If $\chi(F_2)=0$, then $F_2$ consists of annuli, since
$F_2$ cannot contain a M\"obius band component.
Then $\ell(\sigma)=2$.

We regard $F_2$ as the union of $F\cap U_2$ and the rectangles on $A$ again.
Let $E$ be an annulus component of $F_2$, and let $\partial E=e_1\cup e_2$.
Let $g$ be a loop on $\partial U_2$, which is the union of two arcs of $F\cap A$ and two arcs of $\partial F\cap A_2$
and which meets both components of $\partial E$.
Then either $g$ bounds a disk in $E$, or $g$ is essential on $E$.
In the latter case, we replace $g\cap e_1$ with $e_1-\mathrm{Int}(g\cap e_1)$.
Then the resulting $g$ bounds a disk in $E$.
Since $g$ meets the core of $A_2$ in two point (with the same sign), $g$ must be a meridian of $U_2$.
But this implies $p=2$, a contradiction.
Thus we have shown that $\chi(F_2)<0$.
\end{proof}

\begin{proof}[Proof of Theorem \ref{main}(1)]
By Proposition \ref{reduction-even}, $r=pq$.
Consider $r$-Dehn surgery $K(r)$ on $K$.
Then $K(r)=L(p,q)\sharp L(q,p)$ by \cite{M}.
We remark that $L(q,p)$ cannot contain a closed non-orientable surface since $H_1(L(q,p))$ is odd torsion \cite{BW}.
Let $\widehat{F}$ be the closed non-orientable surface obtained by capping $\partial F$ off by a meridian
disk of the attached solid torus of $K(r)$.
Then $L(p,q)$ contains a closed non-orientable surface whose genus is $\widehat{F}$ by \cite[Corollary 5.2]{BW}.
Hence $N(p,q)\le c(K)$.

On the other hand, $N(p,q)$ is also the minimal genus of non-orientable surfaces
in a solid torus having a single $(p,q)$ loop as boundary \cite{BW}.
This means that $K$ bounds a non-orientable surface with genus $N(p,q)$, and thus $c(K)\le N(p,q)$.
Therefore we have $c(K)=N(p,q)$.
This completes the proof of Theorem \ref{main}(1).
\end{proof}

\section{The case where $K$ is odd}\label{odd}

In this section, we prove Theorem \ref{main}(2).

Suppose that $K$ is odd.
Recall that $p>q$.
We first determine which of $N(pq-1,p^2)$ and $N(pq+1,p^2)$ is smaller.
For it, recall the definition of $N(x,y)$ \cite{BW}.
(Although only the case where $x$ is even is considered in \cite{BW}, $N(x,y)$ can be defined in general.)

Let $x$ and $y$ be coprime positive integers.
We write $x/y$ as a continued fraction:

$$\frac{x}{y}=[a_0,a_1,a_2,\dots,a_n]\\
      =a_0+\cfrac{1}{a_1
          +\cfrac{1}{a_2
          +\cfrac{1}{\ddots
          +\cfrac{1}{a_n}}}}$$
where the $a_i$ are integers, $a_0\ge 0$, $a_i>0$ for $1\le i\le n$, and $a_n>1$.
Note that such an expression is unique (cf. \cite{HW}).
Then we define $b_i$ inductively as follows:

\[
b_0=a_0,
\]

\[
b_i=
\begin{cases}
a_i & \text{if $b_{i-1}\ne a_{i-1}$ or if $\displaystyle\sum_{j=0}^{i-1}b_j$ is odd}, \\
0   & \text{if $b_{i-1}=a_{i-1}$ and $\displaystyle\sum_{j=0}^{i-1}b_j$ is even}.
\end{cases}
\]

Let $x_i/y_i=[a_0,a_1,a_2,\dots,a_i]$ be the $i$-th convergent of $x/y$ for $0\le i\le n$.
Thus $x/y=x_n/y_n$.
Let
$$\sum\left(\frac{x_i}{y_i}\right)=\sum_{j=0}^i b_j,$$
and let

$$N(x,y)=\frac{1}{2}\sum\left(\frac{x_n}{y_n}\right)=\frac{1}{2}\sum\left(\frac{x}{y}\right).$$

That is, $\sum(x/y)$ is obtained by adding the $a_i$
successively except when a partial sum is even we skip the next $a_i$.
\cite[Theorem 6.1]{BW} showed that $N(x,y)$ gives the minimal genus of closed non-orientable
surfaces contained in the lens space $L(x,y)$.
(In this case, $x$ must be even. Otherwise, the lens space cannot contain a closed non-orientable surface.
Thus $\sum(x/y)$ must be even when $x$ is even.)

\begin{lemma}\label{recipro}
Let $q/p=[a_0,a_1,a_2,\dots,a_n]$, and let $b_i$ $(0\le i\le n)$ be defined as above.
Then
\[
(pq-1)/p^2=
\begin{cases}
        {[a_0,a_1,a_2,\dots,a_{n-1},a_{n}+1,a_{n}-1,a_{n-1},a_{n-2},\dots,a_2,a_1]} & \text{if $n$ is odd},\\
        {[a_0,a_1,a_2,\dots,a_{n-1},a_{n}-1,a_{n}+1,a_{n-1},a_{n-2},\dots,a_2,a_1]} & \text{if $n$ is even},
\end{cases}
\]

and

\[
(pq+1)/p^2=
\begin{cases}
        {[a_0,a_1,a_2,\dots,a_{n-1},a_{n}-1,a_{n}+1,a_{n-1},a_{n-2},\dots,a_2,a_1]} & \text{if $n$ is odd},\\
        {[a_0,a_1,a_2,\dots,a_{n-1},a_{n}+1,a_{n}-1,a_{n-1},a_{n-2},\dots,a_2,a_1]} & \text{if $n$ is even}.
\end{cases}
\]

\end{lemma}
\begin{remark}
In Lemma \ref{recipro}, $a_0=0$ since $p>q$.
Also, if $a_1=1$, then $[a_0,\dots,a_3,a_2,a_1]$ should be considered to be $[a_0,\dots,a_3,a_2+1]$.
\end{remark}

\begin{proof}
Since $a_0=0$, $p/q=[a_1,a_2,\dots,a_n]$.
We show that $$[a_1,a_2,\dots,a_{n-1},a_{n}+1,a_{n}-1,a_{n-1},a_{n-2},\dots,a_2,a_1]=\frac{p^2}{pq+(-1)^n}.$$
Let $p_i/q_i=[a_1,\dots,a_i]$ be the $i$-th convergent of $[a_1,a_2,\dots,a_n]$ $(1\le i\le n)$.
Then $p_n=p, q_n=q$, and $p_i=a_ip_{i-1}+p_{i-2}$ $(i\ge 3)$.
By induction, $p_{n-1}/p_{n-2}=[a_{n-1},a_{n-2},\dots,a_2,a_1]$.
Then
\begin{align*}
[a_{n}+1,a_{n}-1,a_{n-1},a_{n-2},\dots,a_2,a_1]&=[a_{n}+1,a_{n}-1,[a_{n-1},a_{n-2},\dots,a_2,a_1]]\\
                                               &=[a_{n}+1,a_{n}-1,p_{n-1}/p_{n-2}]\\
                                               &=a_{n}+1+\frac{p_{n-1}}{a_{n}p_{n-1}-p_{n-1}+p_{n-2}}\\
                                               &=a_n+1+\frac{p_{n-1}}{p_n-p_{n-1}}\\
                                               &=\frac{a_np_n-a_np_{n-1}+p_n}{p_n-p_{n-1}}\\
                                               &=\frac{a_np_n+p_{n-2}}{p_n-p_{n-1}}.
\end{align*}
Let us denote this by $c_1$.
Next, let $c_2=[a_{n-1},c_1]$.
Then 
$$c_2=\frac{(a_{n-1}a_{n}+1)p_n-p_{n-3}}{a_np_n+p_{n-2}}.$$
Inductively, we define $c_i=[a_{n-i+1},c_{i-1}]$.
For example, 
$$c_3=[a_{n-2},c_2]=\frac{(a_{n-2}(a_{n-1}a_n+1)+a_n)p_n+p_{n-4}}{(a_{n-1}a_n+1)p_n-p_{n-3}},$$
$$c_4=[a_{n-3},c_3]=\frac{(a_{n-3}(a_{n-2}(a_{n-1}a_n+1)+a_n)+(a_{n-1}a_n+1))p_n-p_{n-5}}{(a_{n-2}(a_{n-1}a_n+1)+a_n)p_n+p_{n-4}}.$$
We will show $c_n=p^2/(pq+(-1)^n)$.

Let $p_i'/q_i'=[a_n,a_{n-1},\dots,a_{n-i+1}]$ be the $i$-th convergent of $[a_n,a_{n-1},\dots,a_2,a_1]$, $1\le i\le n$.
Then $p_n'=p_n=p, q_n'=p_{n-1}$, since $p_n/p_{n-1}=[a_n,a_{n-1},\dots,a_2,a_1]$.
Furthermore, 
\begin{align*}
p_1'=a_n, && p_2'=a_{n-1}a_n+1, && p_i'=a_{n-i+1}p_{i-1}'+p_{i-2}',\\
q_1'=1,   && q_2'=a_{n-1},     && q_i'=a_{n-i+1}q_{i-1}'+q_{i-2}'.
\end{align*}
Then 
$$c_1=\frac{p_1'p_n+p_{n-2}}{p_n-p_{n-1}},\quad c_2=\frac{p_2'p_n-p_{n-3}}{p_1'p_n+p_{n-2}},\dots,c_{n-2}=\frac{p_{n-2}'p_n+(-1)^{n-1}p_1}{p_{n-3}'p_n+(-1)^{n}p_2}.$$
Thus 
$$c_{n-1}=[a_2,c_{n-2}]=\frac{p_{n-1}'p_n+(-1)^{n}}{p_{n-2}'p_n+(-1)^{n-1}p_1},$$
since $p_2=a_1a_2+1$.
Finally, 
$$c_n=[a_1,c_{n-1}]=\frac{p_n'p_n}{p_{n-1}'p_n+(-1)^{n}}.$$
Recall that $p_n'=p_n=p$.
Also we can show that $p_n'/p_{n-1}'=[a_1,a_2,\dots,a_n]=p/q$ inductively, and so $p_{n-1}'=q$.
Therefore
$c_n=p^2/(pq+(-1)^n)$
as desired.

Similarly, we can show that
$$[a_1,a_2,\dots,a_{n-1},a_{n}-1,a_{n}+1,a_{n-1},a_{n-2},\dots,a_2,a_1]=\frac{p^2}{pq+(-1)^{n-1}}.$$
This completes the proof of Lemma \ref{recipro}.
\end{proof}

Consider the equation $xq\equiv -1 \pmod p$.
Since $p$ and $q$ are coprime, this has a solution.
In general, if $x$ is a solution, then so is $x+p$.  Hence the parity of the solution $x$
is not well-defined, because $p$ is odd.
But the equation has the unique solution $x$ satisfying $0<x<p$.
(For, if $x$ and $y$ are such solutions, then $(x-y)q\equiv 0 \pmod p$.
Then $x \equiv y \pmod p$, and so $x=y$.)
If such $x$ is even (resp.~odd), then the ordered pair $(p,q)$
is said to be of \textit{type\/} $A$ (resp.~$B$).

Let $q/p=[a_0,a_1,a_2,\dots,a_n]$, and let $q_i/p_i=[a_0,a_1,a_2,\dots,a_i]$ be the $i$-th convergent for $0\le i\le n$.
We will use $b_i$, $0\le i\le n$, defined for $q/p=[a_0,a_1,a_2,\dots,a_n]$ as before.

\begin{lemma}\label{criterion}
If $n$ is even \textup{(}resp.~odd\textup{)},
then the pair $(p,q)$ is of type A if and only if $p_{n-1}$ is even \textup{(}resp.~odd\textup{)}.
\end{lemma}

\begin{proof}
It is well known that $q_np_{n-1}-q_{n-1}p_n=(-1)^{n-1}$.
Recall $p_n=p, q_n=q$.
If $n$ is even, then $p_{n-1}q\equiv -1 \pmod{p}$.
Since $0<p_{n-1}<p_n=p$ (\cite{HW}), $p_{n-1}$ gives the unique solution of $xq\equiv -1 \pmod p$ such that $0<x<p$.

If $n$ is odd, then $-p_{n-1}q\equiv -1 \pmod{p}$.
Hence $p-p_{n-1}$ gives the unique solution of $xq\equiv -1 \pmod p$ such that $0<x<p$.
Since $p$ is odd, $p_{n-1}$ and $p-p_{n-1}$ have distinct parities.
\end{proof}

\begin{remark}
In Lemma \ref{criterion}, $p_{n-1}$ and $q_{n-1}$ have distinct parities, and hence
when $n$ is odd, $(p,q)$ is of type A if and only if $q_{n-1}$ is even.
\end{remark}

\begin{lemma}\label{odd-pn-1}
If $p_{n-1}$ is odd, then $b_{n-1}=a_{n-1}$ and $\sum(q_{n-1}/p_{n-1})$ is even.
\end{lemma}

\begin{proof}
Since $p_{n-1}$ and $q_{n-1}$ have distinct parities, $q_{n-1}$ is even.
Let us consider $q_{n-1}/p_{n-1}=[a_0,a_1,\dots,a_{n-1}]$.
(If $a_{n-1}=1$, this is $[a_0,a_1,\dots,a_{n-2}+1]$.)
By \cite{BW}, $\sum(q_{n-1}/p_{n-1})$ is even, because $q_{n-1}$ is even.
(In fact, $\sum(q_{n-1}/p_{n-1})$ is twice the minimal genus
of closed non-orientable surfaces contained in the lens space $L(q_{n-1},p_{n-1})$.)

Now, recall another interpretation of $\sum(q_{n-1}/p_{n-1})$ \cite[p107]{BW}.
Consider the ``step'' changing $[a_0,a_1,\dots,a_{n-1}]$ ($=[a_0,a_1,\dots,a_{n-2}+a_{n-1}]=[a_0,a_1,\dots,a_{n-2}+1]$ if $a_{n-1}=1$) as follows:

\[
[a_0,a_1,\dots,a_{n-1}]\to
\begin{cases}
         {[a_0,a_1,\dots,a_{n-1}-2]} & \text{if $a_{n-1}\ge 4$},\\
         {[a_0,a_1,\dots,a_{n-2},1]=[a_0,a_1,\dots,a_{n-2}+1]} & \text{if $a_{n-1}=3$},\\
         {[a_0,a_1,\dots,a_{n-3}]} & \text{if $a_{n-1}=2$}.
\end{cases}
\]

Then $\sum(q_{n-1}/p_{n-1})$ is twice the number of steps required to reduce
$[a_0,a_1,\dots,a_{n-1}]$ to $[0]$.
(Such a reduction works only when $q_{n-1}$ is even.)
In other words, the counting process defining  $\sum(q_{n-1}/p_{n-1})$
can be done from either end of the sequence $a_0,a_1,\dots,a_{n-1}$.
This implies that the calculation of $\sum(q_{n-1}/p_{n-1})$ involves the last term $a_{n-1}$.
Hence $b_{n-1}=a_{n-1}$.
\end{proof}

\begin{lemma}\label{oddtype}
If $p_{n-1}$ is odd, then $N(pq-1,p^2)=N(pq+1,p^2)+(-1)^n$.
\end{lemma}

\begin{proof}
Assume $n$ is odd.
Then 
\begin{align*}
(pq-1)/p^2 &=[a_0,a_1,a_2,\dots,a_{n-1},a_{n}+1,a_{n}-1,a_{n-1},a_{n-2},\dots,a_2,a_1],\\
(pq+1)/p^2 &=[a_0,a_1,a_2,\dots,a_{n-1},a_{n}-1,a_{n}+1,a_{n-1},a_{n-2},\dots,a_2,a_1]
\end{align*}
by Lemma \ref{recipro}.
In the calculations of $\sum((pq-1)/p^2)$ and $\sum((pq+1)/p^2)$,
$b_{n-1}=a_{n-1}$ and the partial sum $\sum_{j=0}^{n-1}b_j$ is even by Lemma \ref{odd-pn-1}.
Hence the next term is skipped, that is, $b_n=0$.
Thus we see that $\sum((pq-1)/p^2)=\sum((pq+1)/p^2)-2$.
Then $N(pq-1,p^2)=N(pq+1,p^2)-1$.

Assume $n$ is even.
Then a similar argument shows that $N(pq-1,p^2)=N(pq+1,p^2)+1$.
\end{proof}

\begin{lemma}\label{eventype}
If $p_{n-1}$ is even, then $N(pq-1,p^2)+(-1)^n=N(pq+1,p^2)$.
\end{lemma}

\begin{proof}
Since $q_{n-1}/p_{n-1}=[a_0,a_1,\dots,a_{n-1}]=[0,a_1,\dots,a_{n-1}]=1/[a_1,a_2,\dots,a_{n-1}]$, we see
$p_{n-1}/q_{n-1}=[a_1,a_2,\dots,a_{n-1}]$.
Thus $\sum(p_{n-1}/q_{n-1})$ is even and
its calculation involves the last term $a_{n-1}$, since $p_{n-1}$ is even, 
as in the proof of Lemma \ref{odd-pn-1}.

Assume $n$ is odd.
Recall 
\begin{align*}
(pq-1)/p^2 &=[a_0,a_1,a_2,\dots,a_{n-1},a_{n}+1,a_{n}-1,a_{n-1},a_{n-2},\dots,a_2,a_1],\\
(pq+1)/p^2 &=[a_0,a_1,a_2,\dots,a_{n-1},a_{n}-1,a_{n}+1,a_{n-1},a_{n-2},\dots,a_2,a_1].
\end{align*}
For the calculations of $\sum((pq-1)/p^2)$ and $\sum((pq+1)/p^2)$,
we use the counting procedure backward as stated in the proof of Lemma \ref{odd-pn-1}.
Starting from $a_1$, it reaches $[a_0,a_1,\dots,a_{n-1},a_{n}\pm 1,a_{n}\mp 1,0]=[a_0,a_1,\dots,a_{n-1},a_n\pm 1]$.
Thus we see that $\sum((pq-1)/p^2)-2=\sum((pq+1)/p^2)$.
Then $N(pq-1,p^2)-1=N(pq+1,p^2)$.

When $n$ is even, a similar argument shows $N(pq-1,p^2)+1=N(pq+1,p^2)$.
\end{proof}

We are ready to determine which of $N(pq-1,p^2)$ and $N(pq+1,p^2)$ is smaller.

\begin{proposition}\label{final}
\begin{itemize}
\item[(1)] If the pair $(p,q)$ is of type A, then $N(pq-1,p^2)+1=N(pq+1,p^2)$.
\item[(2)] If the pair $(p,q)$ is of type B, then $N(pq+1,p^2)+1=N(pq-1,p^2)$.
\end{itemize}
\end{proposition}

\begin{proof}
Suppose $(p,q)$ is of type A.
If $n$ is even, then so is $p_{n-1}$ by Lemma \ref{criterion}.
Then $N(pq-1,p^2)+1=N(pq+1,p^2)$ by Lemma \ref{eventype}.
If $n$ is odd, then $p_{n-1}$ is odd.
Then $N(pq-1,p^2)=N(pq+1,p^2)-1$ by Lemma \ref{oddtype}.
This proves (1).
(2) can be proved similarly.
\end{proof}

\begin{remark}
In general, the difference between $N(x,y)$ and $N(x+2,y)$ can be big.
For example, $N(26,25)=13, N(28,25)=6$.
\end{remark}

\begin{figure}[tb]
\includegraphics*[scale=0.45]{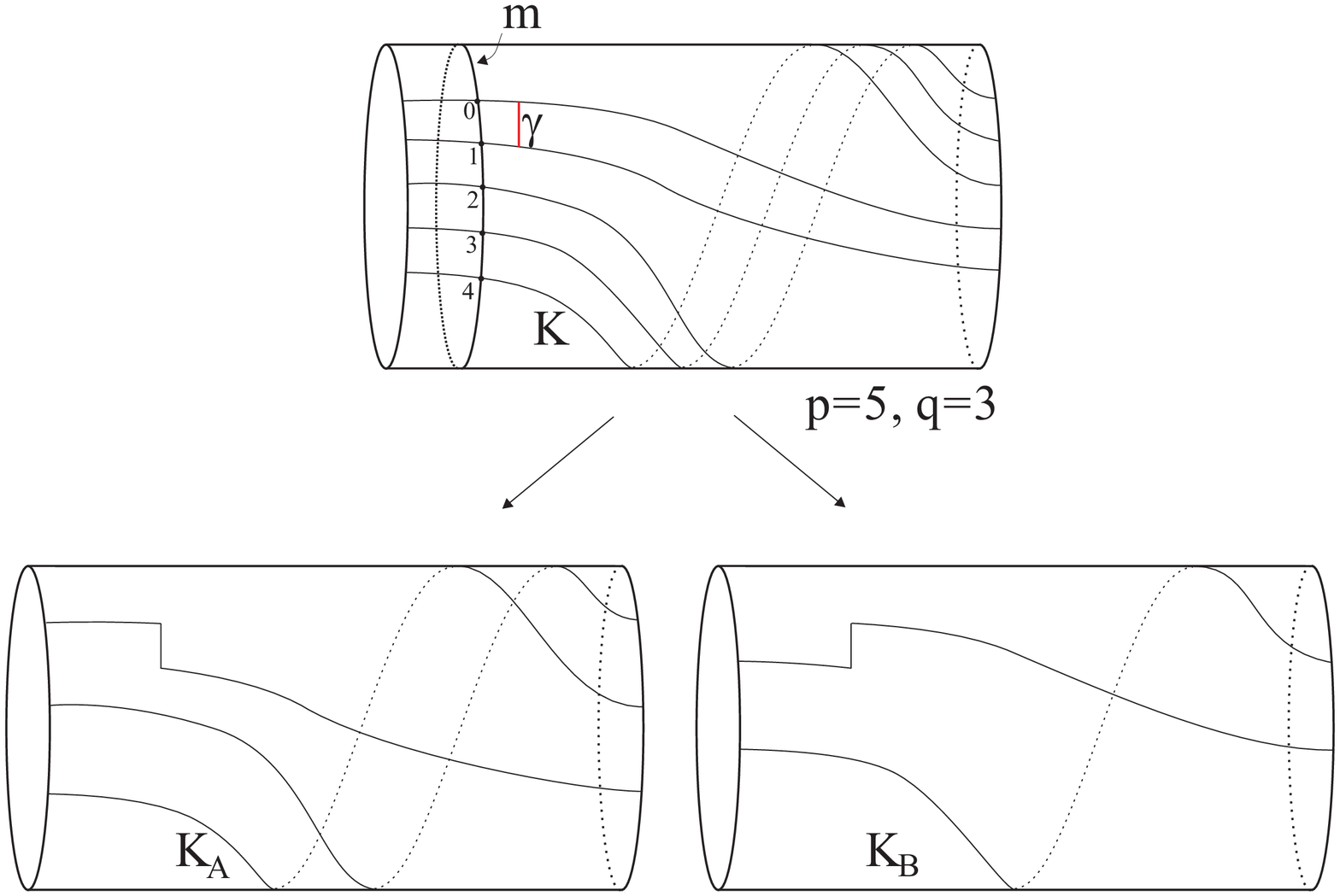}
\caption{}\label{gamma}
\end{figure}

Now, we will explain the geometric meaning of types A and B.
Recall that $K=T(p,q)$ lies on the standard torus $T$, which decomposes $S^3$ into
two solid tori $V$ and $W$, and that
$K$ runs $p$ times longitudinally and $q$ times meridionally with respect to $V$.
Choose an arc $\gamma$ on $T$ as shown in Figure \ref{gamma}.
(Here, the end circles of the cylinder are identified to form $T$.)
Then $\partial \gamma$ splits $K$ into two arcs $K_a$ and $K_b$.
Let $K_A=\gamma \cup K_a$ and $K_B=\gamma\cup K_b$ as shown in Figure \ref{gamma}.
Clearly, both $K_A$ and $K_B$ are torus knots, which are uniquely determined by $K$.
In Figure \ref{gamma}, $K=T(5,3)$, $K_A=T(3,2)$ and $K_B=T(2,1)$ (with respect to $V$).

Let $K_A=T(r_1,s_1)$ and $K_B=T(r_2,s_2)$.
Then $p=r_1+r_2$ and $q=s_1+s_2$.

\begin{lemma}\label{mq}
\begin{itemize}
\item[(1)] If $(p,q)$ is of type $A$, then $r_1$ and $s_2$ are even, and hence $s_1$ and $r_2$ are odd.
\item[(2)] If $(p,q)$ is of type $B$, then $s_1$ and $r_2$ are even, and hence $r_1$ and $s_2$ are odd.
\end{itemize}
\end{lemma}

\begin{proof}
Let $m$ be a meridian of $V$ near $\gamma$, and
label the points of $K\cap m$, $0,1,2,\dots,p-1$ along $m$ as shown in Figure \ref{gamma}.
Start the point with label $0$ and follow $K_A$ in the direction of $\gamma$.
Then we will come back to the point with label $1+q$ after running once longitudinally.
Hence $K_A$ gives the equation $1+r_1q\equiv 0\pmod p$.
Thus if $(p,q)$ is of type $A$, then $r_1$ is even, and $s_1$ is odd.
Since $p=r_1+r_2$ and $q=s_1+s_2$ are odd, $r_2$ is odd and $s_2$ is even.
This proves (1). (2) follows similarly.
\end{proof}

\begin{lemma}\label{tunnel}
\begin{itemize}
\item[(1)] If $(p,q)$ is of type $A$, then $K$ bounds a non-orientable surface
with genus $N(r_1,s_1)+N(s_2,r_2)$ and boundary slope $pq-1$.
\item[(2)] If $(p,q)$ is of type $B$, then $K$ bounds a non-orientable surface
with genus $N(s_1,r_1)+N(r_2,s_2)$ and boundary slope $pq+1$.
\end{itemize}
\end{lemma}

\begin{proof}
Assume that $(p,q)$ is of type $A$.
By Lemma \ref{mq}, $K_A$ bounds a non-orientable surface $F_A$ contained in $V$ with genus $N(r_1,s_1)$,
and $K_B$ bounds such $F_B$ in $W$ with genus $N(s_2,r_2)$.
Let $G=F_A\cup F_B$.
Since $F_A\cap F_B=\gamma$, $\partial G=K$ and $G$ has genus $N(r_1,s_1)+N(s_2,r_2)$.
It is easy to see that $G$ has the desired boundary slope from Figure \ref{AB}.
This proves (1). 

\begin{figure}[tb]
\includegraphics*[scale=0.45]{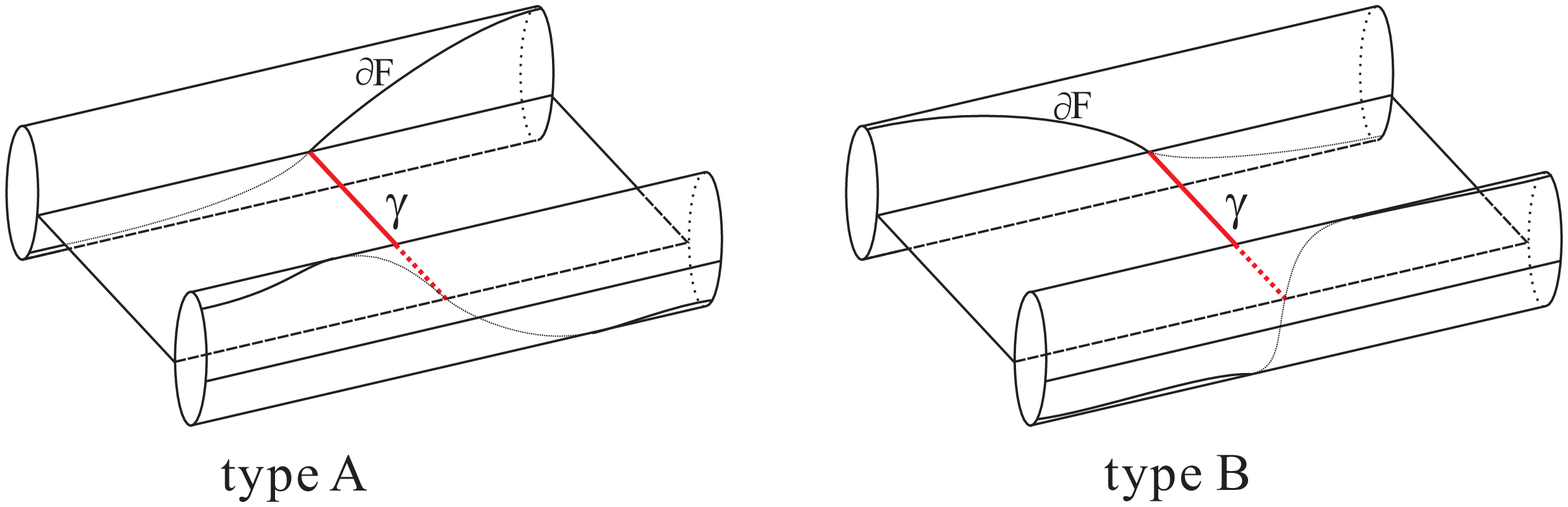}
\caption{}\label{AB}
\end{figure}

A similar argument shows (2).
\end{proof}

Recall that $q_{n-1}/p_{n-1}=[a_0,a_1,\dots,a_{n-1}]$ is
the $(n-1)$ th convergent of $q_n/p_n=q/p=[a_0,a_1,\dots,a_n]$.

\begin{lemma}\label{KAKB}
If $n$ is even, then
$r_1=p_{n-1}$ and $s_1=q_{n-1}$.
If $n$ is odd, $r_2=p_{n-1}$ and $s_2=q_{n-1}$.
\end{lemma}

\begin{proof}
Suppose $n$ is even.
As shown in the proof of Lemma \ref{criterion}, we see $r_1=p_{n-1}$.
Recall that $q_{n-1}p-p_{n-1}q=1$ and $0<q_{n-1}<q_n=q$.
By considering the intersection number between $K$ and $K_A$ on $T$ (after
giving $K$ an orientation, which induces that of $K_A$),
we have $|s_1p-p_{n-1}q|=1$.
Since $p\ge 5$, we have $(q_{n-1}-s_1)p=0$, and then $s_1=q_{n-1}$.

When $n$ is odd, we see $r_1=p-p_{n-1}$.
Then a similar argument shows that $s_1=q-q_{n-1}$.
\end{proof}

\begin{lemma}\label{claim}
$[a_1,a_2,\dots,a_{n-1},a_n-1]=(p-p_{n-1})/(q-q_{n-1})$.
\end{lemma}

\begin{proof}
Recall that $p_{n-1}/q_{n-1}=[a_1,a_2,\dots,a_{n-1}]$.
Put $P/Q=[a_1,a_2,\dots,a_{n-1},a_n-1]$.
If $a_n>2$, then
$P=(a_n-1)p_{n-1}+p_{n-2}=p-p_{n-1}$ as desired.
If $a_n=2$, then $P/Q=[a_1,a_2,\dots,a_{n-1}+1]$.
Thus $P=(a_{n-1}+1)p_{n-2}+p_{n-3}=p_{n-1}+p_{n-2}$.
Since $p=p_n=a_np_{n-1}+p_{n-2}=2p_{n-1}+p_{n-2}$,
$p_{n-1}+p_{n-2}=p-p_{n-1}$ as desired again.
Similarly for $Q$.
\end{proof}

\begin{lemma}\label{sum}
\begin{itemize}
\item[(1)] If $(p,q)$ is of type $A$, then $N(pq-1,p^2)=N(r_1,s_1)+N(s_2,r_2)$.
\item[(2)] If $(p,q)$ is of type $B$, then $N(pq+1,p^2)=N(s_1,r_1)+N(r_2,s_2)$.
\end{itemize}
\end{lemma}

\begin{proof}
(1) Suppose that $p_{n-1}$ is odd.
By Lemma \ref{criterion}, $n$ is odd.
Then $K_B=T(p_{n-1},q_{n-1})$ by Lemma \ref{KAKB}, and hence
$K_A=T(p-p_{n-1},q-q_{n-1})$.
Recall that 
$(pq-1)/p^2=[a_0,a_1,a_2,\dots,a_{n-1},a_{n}+1,a_{n}-1,a_{n-1},a_{n-2},\dots,a_2,a_1]$
by Lemma \ref{recipro}.
In the counting procedure of $\sum((pq-1)/p^2)$,
$\sum_{j=0}^{n-1}b_j=\sum(q_{n-1}/p_{n-1})$ by Lemma \ref{odd-pn-1}.

As stated in the proof of Lemma \ref{odd-pn-1},
the counting procedure can be done backward.
But $[a_1,a_2,\dots,a_{n-1},a_n-1]=(p-p_{n-1})/(q-q_{n-1})$ by Lemma \ref{claim}.
Thus $$\sum((pq-1)/p^2)=\sum(q_{n-1}/p_{n-1})+\sum((p-p_{n-1})/(q-q_{n-1})).$$
Hence we have $N(pq-1,p^2)=N(s_2,r_2)+N(r_1,s_1)$.

Suppose that $p_{n-1}$ is even.
Then $n$ is even by Lemma \ref{criterion} so that $K_A=T(p_{n-1},q_{n-1})$.
Recall
$(pq-1)/p^2=[a_0,a_1,a_2,\dots,a_{n-1},a_{n}-1,a_{n}+1,a_{n-1},a_{n-2},\dots,a_2,a_1]$.
Perform the counting of $\sum((pq-1)/p^2)$ backward.
Then $\sum[a_1,a_2,\dots,a_{n-1}]=\sum(p_{n-1}/q_{n-1})$, and
$\sum[a_0,a_1,\dots,a_{n-1},a_n-1]=\sum((q-q_{n-1})/(p-p_{n-1}))$.
Hence $N(pq-1,p^2)=N(s_2,r_2)+N(r_1,s_1)$ again.

A similar argument shows (2).
\end{proof}

\begin{corollary}\label{half}
$c(K)\le N(pq-1,p^2)$ \textup{(}resp.~$N(pq+1,p^2)$\textup{)}
when $K$ is of type A \textup{(}resp.~B\textup{)}.
\end{corollary}

\begin{proof}
By Lemmas \ref{tunnel} and \ref{sum}, if $(p,q)$ is of type A (resp.~B), 
then $K$ can bound a non-orientable surface with genus $N(pq-1,p^2)$
(resp.~$N(pq+1,p^2)$).
Hence $c(K)\le N(pq-1,p^2)$ (resp.~$N(pq+1,p^2)$) when $K$ is of type A (resp.~B).
\end{proof}

We now determine the boundary slope $r$ of a minimal genus non-orientable spanning
surface $F$ of $K$.

\begin{proposition}\label{reduction-odd}
$r=pq\pm 1$.
\end{proposition}

\begin{proof}
Suppose not.
Since $\Delta=|r-pq|$ is odd, $\Delta\ge 3$.
The proof is divided into two cases.

\textit{Case 1. $\Delta\ge 5$.}

We perform $(\Delta-1)/2$ disk splittings for $F$ along mutually disjoint $(\Delta-1)/2$ rectangles on $A$
as in the proof of Proposition \ref{reduction-even}.
By this operation, $\partial F$ breaks up into a loop with slope $pq\pm 1$ and 
the others inessential on $\partial E(K)$.
It can be observed that the latter is not empty under the condition $\Delta\ge 5$.
Let $F_1$ be the resulting connected surface whose boundary contains a component of slope $pq\pm 1$,
and let $F_2$ be the others.
Note that $\chi(F)=\chi(F_1)+\chi(F_2)$.
Thus $F_1$ is non-orientable, and $F_2$ may be empty or disconnected.
Also $F_1\cap A\ne \emptyset$, but $F_2\cap A=\emptyset$.

Now, $A$ splits $E(K)$ into two solid tori $U_1, U_2$, where we assume
$F_2\subset U_2$.
Let $A_i=\partial E(K)\cap U_i$.
Define $\sigma$ by using the arcs $\partial F\cap A_2$ as in the proof of Proposition \ref{reduction-even}.

\begin{claim}\label{single}
$F_2\ne \emptyset$ and $F_1$ has a single boundary component.
\end{claim}

\begin{proof}[Proof of Claim \ref{single}]
The argument in the proof of Claim \ref{F2} works, and so $F_2$ has no disk components.
Thus $\chi(F)=\chi(F_1)+\chi(F_2)\le \chi(F_1)$.
If either $F_2=\emptyset$ or $|\partial F_1|>1$, then $\partial F_1$ contains an inessential component on $\partial E(K)$.
Then cap all inessential components of $\partial F_1$ off by disks on $\partial E(K)$.
The resulting surface gives a non-orientable spanning surface of $K$
with fewer betti number than $F$.  This contradicts the minimality of $F$.
\end{proof}

Recall that $\ell(\sigma)$ is the (common) length of the orbits of $\sigma$.

\begin{claim}\label{F1}
$\ell(\sigma)=1$.
\end{claim}

\begin{proof}[Proof of Claim \ref{F1}]
Suppose not.
Let $\alpha$ be the arc of $F\cap A$ which is disjoint from all rectangles used in the
disk splittings.
This arc remains in $F_1$.
We choose the labeling so that the point $\alpha\cap a_1$ has $\Delta$.
Since $\ell(\sigma)\ne 1$, the point $\alpha\cap a_1$
is not connected with $\alpha\cap a_2$ (with label $\Delta$) by an arc in $A_2$.
See Figure \ref{orbit1} where $\Delta=7$.
(In this and successive figures, we assume $r=pq+\Delta$.  But
the situation when $r=pq-\Delta$ is similar.)

\begin{figure}[tb]
\includegraphics*[scale=0.45]{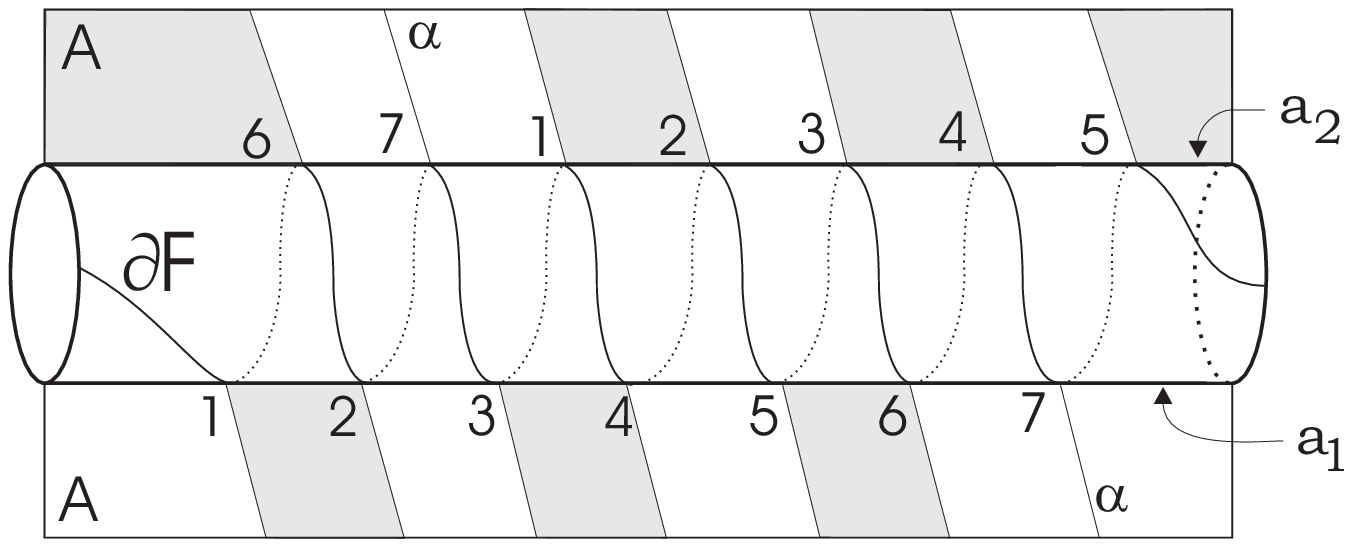}
\caption{}\label{orbit1}
\end{figure}

If either situation as shown in Figure \ref{orbit2} happens, then
the disk splittings give rise to inessential components on both $A_1$ and $A_2$.
Then $\partial F_1$ would have an inessential component, which is impossible by Claim \ref{single}. 
Thus the only possible configurations are as shown in Figure \ref{orbit3}.
But it is easy to see that $\partial F_1$ has an inessential component (and $F_2=\emptyset$) in either case.
This contradicts Claim \ref{single} again.
\end{proof}

\begin{figure}[tb]
\includegraphics*[scale=0.45]{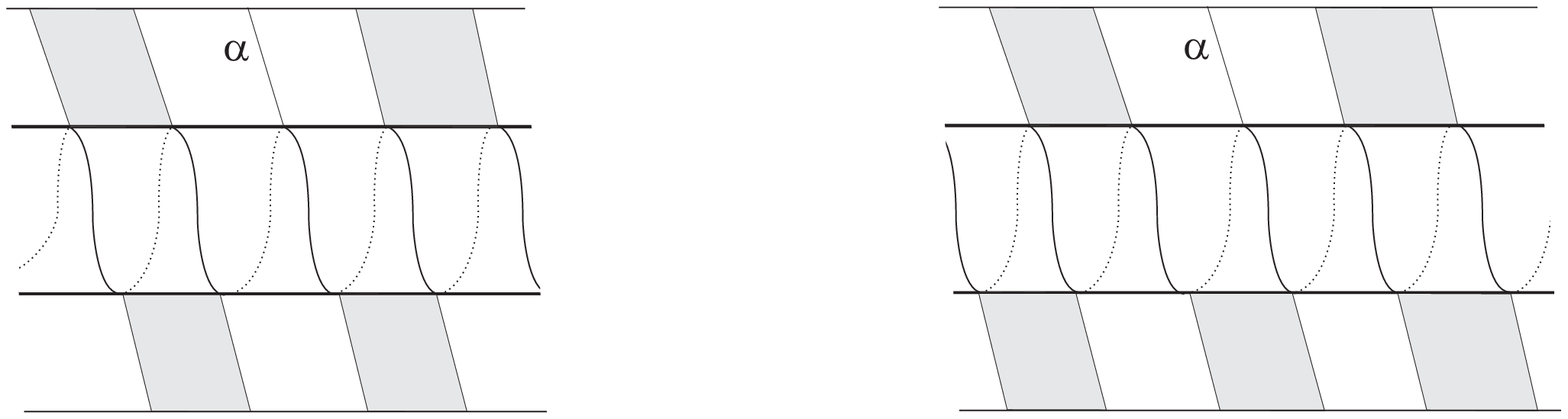}
\caption{}\label{orbit2}
\end{figure}

\begin{figure}[tb]
\includegraphics*[scale=0.45]{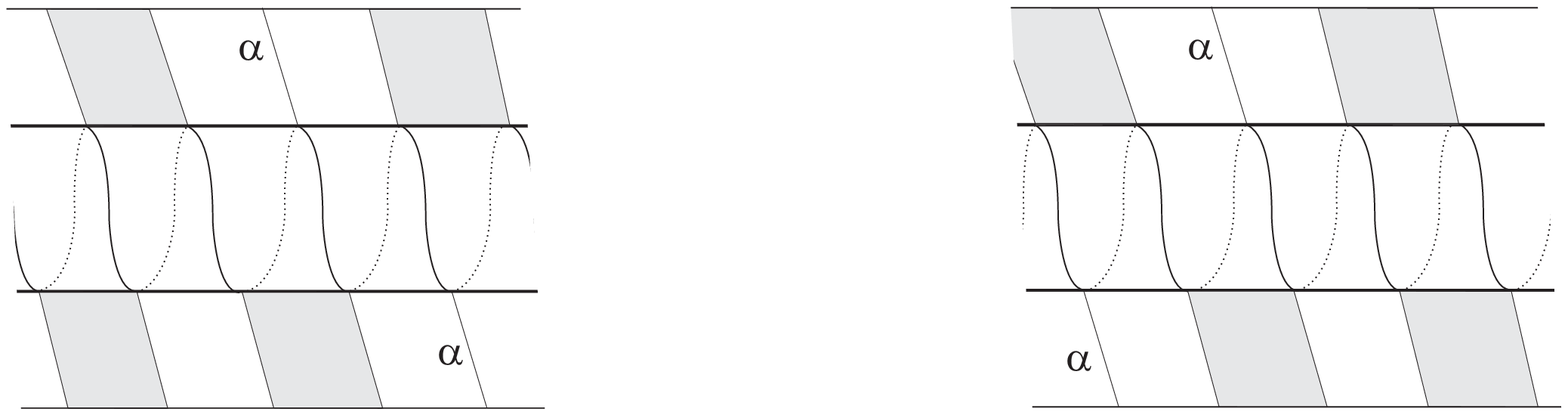}
\caption{}\label{orbit3}
\end{figure}

Thus $|\partial F_2|=(\Delta-1)/2\ge 2$.
If $\chi(F_2)=0$, then $F_2$ consists of annuli.
Then the same argument as in the proof of Claim \ref{F2} gives a contradiction.
(That is, let $f$ be a loop consisting of an arc of $F\cap A$ and one of $\partial F\cap A_2$, which meets a component $E$ of $F_2$.
Then either $f$ bounds a disk in $E$, or $f$ is essential in $E$.
In the latter, replace $f\cap e$ with $e-\mathrm{Int}(f\cap e)$, where $e$ is the component of $\partial E$ meeting $f$.
Then it is necessarily inessential in $E$.)
Thus we have shown that $\chi(F_2)<0$, and so $\chi(F)<\chi(F_1)$, a contradiction.

\textit{Case 2. $\Delta=3$.}

We perform a disk splitting for $F$ along one rectangle on $A$ as before.
Use the same notation as in Case 1.
Then $F_2=\emptyset$ or $|\partial F_2|=1$.
In the latter case, the argument in the proof of Claim \ref{F2} shows that $F_2$ is not a disk.
Thus $\chi(F_2)<0$.   Then $\chi(F)<\chi(F_1)$, a contradiction.
Hence $F_2=\emptyset$, and the only possible configurations are as shown in Figure \ref{orbit4}.

\begin{figure}[tb]
\includegraphics*[scale=0.55]{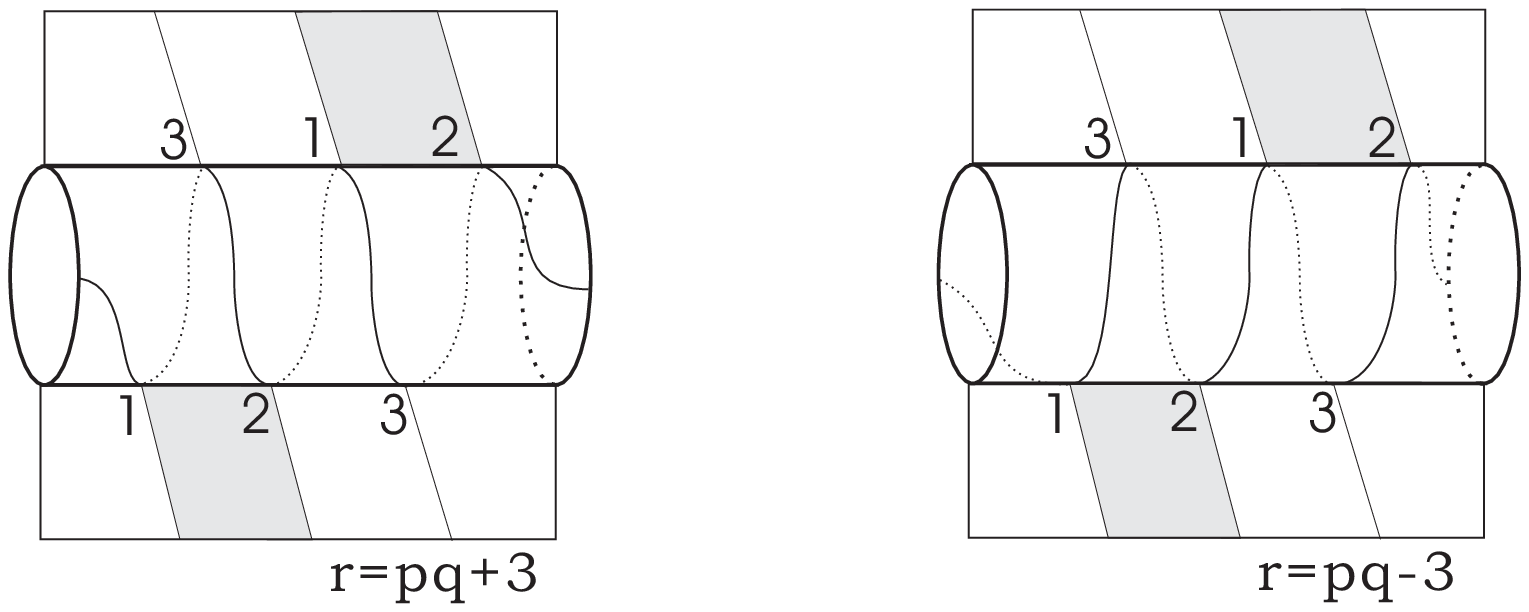}
\caption{}\label{orbit4}
\end{figure}

\begin{claim}\label{onetype}
If $r=pq+3$ \textup{(}resp. $pq-3$\textup{)}, then $(p,q)$ is of type A \textup{(}resp B\textup{)}.
\end{claim}

\begin{proof}[Proof of Claim \ref{onetype}]
Suppose that $r=pq+3$ and that $(p,q)$ is of type B.
Let $F'$ be the resulting surface obtained from $F$ by a disk splitting.
Then $F'$ has boundary slope $pq-1$, and $\chi(F')=\chi(F)$.
Note that $F'$ is connected and non-orientable.
Hence $N(pq-1,p^2)\le c(K)$.
This contradicts that $c(K)\le N(pq+1,p^2)<N(pq-1,p^2)$ (Proposition \ref{final} and Corollary \ref{half}).

The case where $r=pq-3$ is similar.
\end{proof}

We will exclude $r=pq+3$.  The argument to the case $r=pq-3$ is similar.

Recall $E(K)=U_1\cup_A U_2$ and $A_i=\partial E(K)\cap U_i$.
We assume that $K$ has type $(p,q)$ with respect to $U_1$.
Let $G_i=F\cap U_i$ and $K_i=\partial G_i$ for $i=1, 2$.
Then $\chi(G_1)+\chi(G_2)-3=\chi(F)$, and so $\beta_1(G_1)+\beta_1(G_2)+2=\beta_1(F)$.

Let us also recall that $q/p=[a_0,a_1,\dots,a_n]$ and $q_{n-1}/p_{n-1}=[a_0,a_1,\dots,a_{n-1}]$.

\begin{claim}\label{trans}
\begin{itemize}
\item[(1)] If $n$ is odd, $K_1$ has type $(3p_{n-1}+p(3k-1),3q_{n-1}+q(3k-1))$ with respect to $U_1$,
and $K_2$ has type $(-3q_{n-1}+q(2-3k),-3p_{n-1}+p(2-3k))$ with respect to $U_2$, for some even $k$.
\item[(2)] If $n$ is even, $K_1$ has type $(-3p_{n-1}+p(3k-1),-3q_{n-1}+q(3k-1))$ with respect to $U_1$,
and $K_2$ has type $(3q_{n-1}+q(2-3k),3p_{n-1}+p(2-3k))$ with respect to $U_2$, for some odd $k$.
\end{itemize}
\end{claim}

\begin{proof}[Proof of Claim \ref{trans}]
(1) Let $M$ denote the matrix 
\[
\begin{pmatrix}
p_{n-1} & q_{n-1} \\
p  &  q \\
\end{pmatrix}
.\]
Note $\det M=(-1)^{n-1}=1$.
Both $A_i$ are parallel to the annulus $T-\mathrm{Int}A$ in $N(K)$.
Thus $K_i$ is isotopic to $K_i'$, say, on $T$ rel $A$.
Let $f$ be an orientation-preserving self-homeomorphism on $T$ corresponding to $M$.
Then it is easy to see from Figure \ref{orbit4} that $f^{-1}(K_1')$ has type $(3,3k-1)$
and $f^{-1}(K_2')$ has type $(-3,2-3k)$ for some integer $k$ (with suitable orientations).
The conclusion follows easily from this.

Remark that $p_{n-1}$ is odd and $q_{n-1}$ is even by Claim \ref{onetype}, Lemma \ref{criterion} and the remark following it.
Since $K_i$ bounds $G_i$ in $U_i$, we see that $k$ is even.

(2) Use the matrix
\[
\begin{pmatrix}
-p_{n-1} & -q_{n-1} \\
p  &  q \\
\end{pmatrix}
\]
instead of $M$.
\end{proof}

\begin{claim}\label{complicated}
\begin{itemize}
\item[(1)] Let $n$ be odd.  Then 
\[
\frac{3p_{n-1}+p(3k-1)}{3q_{n-1}+q(3k-1)}=
\begin{cases}
       {[a_1,a_2,\dots,a_{n-1},a_{n}-1,1,|k|-1,3]} & \text{if $k\le -2$},\\
       {[a_1,a_2,\dots,a_{n-1},a_{n}-3]} & \text{if $k=0$},\\
       {[a_1,a_2,\dots,a_{n-1},a_{n},k-1,1,2]} & \text{if $k\ge 2$}
\end{cases}
\]

and

\[
\frac{-3q_{n-1}+q(2-3k)}{-3p_{n-1}+p(2-3k)}=
\begin{cases}
     {[0,a_1,a_2,\dots,a_{n-1},a_{n}-1,1,|k|-1,1,2]} & \text{if $k\le -2$},\\
     {[0,a_1,a_2,\dots,a_{n-1},a_{n}-2,2]} & \text{if $k=0$},\\
     {[0,a_1,a_2,\dots,a_{n-1},a_{n},k-1,3]} & \text{if $k\ge 2$.}
\end{cases}
\]
\item[(2)] Let $n$ be even. Then 
\[
\frac{-3p_{n-1}+p(3k-1)}{-3q_{n-1}+q(3k-1)}=
\begin{cases}
    {[a_1,a_2,\dots,a_{n-1},a_{n},|k|,3]} & \text{if $k\le -1$},\\
    {[a_1,a_2,\dots,a_{n-1},a_{n}-2,2]} & \text{if $k=1$},\\
    {[a_1,a_2,\dots,a_{n-1},a_{n}-1,1,k-2,1,2]} & \text{if $k\ge 3$}
\end{cases}
\]

and

\[
\frac{3q_{n-1}+q(2-3k)}{3p_{n-1}+p(2-3k)}=
\begin{cases}
    {[0,a_1,a_2,\dots,a_{n-1},a_{n},|k|,1,2]} & \text{if $k\le -1$},\\
    {[0,a_1,a_2,\dots,a_{n-1},a_{n}-3]} & \text{if $k=1$},\\
    {[0,a_1,a_2,\dots,a_{n-1},a_{n}-1,1,k-2,3]} & \text{if $k\ge 3$.}
\end{cases}
\]
\end{itemize}
\end{claim}

(We consider that $[\dots,a_{n-1},-1]=[\dots,a_{n-1}-1]$,
$[\dots,a_{n-2},a_{n-1},0]=[\dots,a_{n-2}]$,
and $[\dots,a_{n-1},0,2]=[\dots,a_{n-1}+2]$.) 

\begin{proof}[Proof of Claim \ref{complicated}]
This is straightforward.
Suppose that $n$ is odd and $k\le -2$.
Then
\begin{align*}
[a_1,a_2,\dots,a_{n-1},a_n-1,1,|k|-1,3]&=[a_1,a_2,\dots,a_{n-1},[a_n-1,1,|k|-1,3]]\\
   &=\frac{[a_n-1,1,|k|-1,3]p_{n-1}+p_{n-2}}{[a_n-1,1,|k|-1,3]q_{n-1}+q_{n-2}}\\
   &=\frac{(a_n-1)p_{n-1}+p_{n-2}+\frac{3k+2}{3k-1}p_{n-1}}{(a_n-1)q_{n-1}+q_{n-2}+\frac{3k+2}{3k-1}q_{n-1}}\\
   &=\frac{p-p_{n-1}+\frac{3k+2}{3k-1}p_{n-1}}{q-q_{n-1}+\frac{3k+2}{3k-1}q_{n-1}}\\
   &=\frac{(3k-1)p+3p_{n-1}}{(3k-1)q+3q_{n-1}}\\
\end{align*}
as desired.  The other cases can be checked similarly.
\end{proof}

By using these continued fractions, we can evaluate $\beta_1(G_i)$.
Assume $n$ is odd.
As in the proof of Lemma \ref{sum}, 
\[
\beta_1(G_1)\ge N(3p_{n-1}+p(3k-1),3q_{n-1}+q(3k-1))=
\begin{cases}
N(p-p_{n-1},q-q_{n-1})-1 & \text{if $k=0$},\\
N(p-p_{n-1},q-q_{n-1})+\frac{|k|}{2}+1 & \text{otherwise},
\end{cases}
\]

and

\[
\beta_1(G_2)\ge N(-3q_{n-1}+q(2-3k),-3p_{n-1}+p(2-3k))=
\begin{cases}
N(q_{n-1},p_{n-1})+1 & \text{if $k=0$},\\
N(q_{n-1},p_{n-1})+\frac{|k|}{2}+1 & \text{otherwise.}
\end{cases}
\]

Recall that $\beta_1(F)=\beta_1(G_1)+\beta_1(G_2)+2$.
Thus $c(K)=\beta_1(F)>N(pq-1,p^2)$, which contradicts Corollary \ref{half}.

The argument to the case $n$ even is similar.
\end{proof}

\begin{proposition}\label{boundary-slope2}
If $(p,q)$ is of type A \textup{(}resp.~B\textup{)}, then $r=pq-1$
\textup{(}resp.~$pq+1$\textup{)}.
\end{proposition}

\begin{proof}
By Proposition \ref{reduction-odd}, $r=pq-1$ or $pq+1$.
Suppose that $(p,q)$ is of type A.
If $r=pq+1$, then $pq+1$-surgery $K(pq+1)=L(pq+1,p^2)$ contains the closed non-orientable
surface $\widehat{F}$, obtained by capping $\partial F$ off
by a meridian disk of the attached solid torus.
Hence $N(pq+1,p^2)\le c(K)$.
This contradicts Lemma \ref{final} and Corollary \ref{half}.
Thus $r=pq-1$.

A similar argument shows that $r=pq+1$ when $(p,q)$ is of type B.
\end{proof}

\begin{proof}[Proof of Theorem \ref{main}(2)]
Assume $(p,q)$ is of type A.
If $F$ is a non-orientable spanning surface of $K$ realizing $c(K)$,
then its boundary slope $r$ is $pq-1$ by Proposition \ref{boundary-slope2}.
We can cap $\partial F$ off by a meridian disk of the attached solid torus of $K(r)$.
Thus $N(pq-1,p^2)\le c(K)$.
Combined with Corollary \ref{half}, we have $c(K)=N(pq-1,p^2)$.
Similarly, if $(p,q)$ is of type $B$, then $c(K)=N(pq+1,p^2)$.
This completes the proof of Theorem \ref{main}(2).
\end{proof}

\begin{proof}[Proof of Theorem \ref{main2}]
By \cite[Proposition 4.3]{MY}, $c(K_1\sharp K_2)=c(K_1)+c(K_2)$ if and only if
$c(K_1)=\Gamma(K_1)$ and $c(K_2)=\Gamma(K_2)$, where $\Gamma(K_i)=\min \{2g(K_i),c(K_i)\}$.
Let $K=T(p,q)$ be a non-trivial torus knot.
As known well, $g(K)=(p-1)(q-1)/2$.
As in the proof of Lemma \ref{span}, $c(K)\le \min\{(p-1)q/2,(q-1)p/2\}$.
Since it is easy to check that $\min\{(p-1)q/2,(q-1)p/2\}<2g(K)$,
we have $\Gamma(K)=c(K)$.
Thus, if $K_1$ and $K_2$ are torus knots, then $c(K_1\sharp K_2)=c(K_1)+c(K_2)$.

By a standard cut-and-paste argument, $\Gamma(K_1\sharp K_2)=\Gamma(K_1)+\Gamma(K_2)$.
(See also \cite[Section 0]{MY}.)
Thus we have the result inductively.
\end{proof}

\section{Examples and comments}\label{ex}

Let $K$ be a torus knot, and let $F$ be a minimal genus non-orientable surface
bounded by $K$.

\begin{example}
Let $K=T(8,3)$.  Since $K$ is even, the boundary slope of $F$ is $24$.
We have $8/3=[2,1,2]$, so that $b_0=2,b_1=0,b_2=2$.
This gives $\sum(8/3)=4$, and hence $c(K)=N(8,3)=2$.
\end{example}

\begin{example}
Let $K=T(7,5)$. Since $4\cdot 5\equiv -1 \pmod{7}$, the pair $(7,5)$ is of type A.
Thus the boundary slope of $F$ is $34$.
Since $34/49=[0,1,2,3,1,3]$, $c(K)=N(34/49)=(0+2+1+3)/2=3$.

Next, let $K=T(25,9)$. Then $(25,9)$ is of type B.
Then the boundary slope of $F$ is $226$.
Since $226/625=[0,2,1,3,3,1,3,1,2]$, $c(K)=N(226/625)=(0+1+3+1+3+2)/2=5$.
\end{example}

In this paper, we proved that the boundary slope of a minimal genus non-orientable spanning surface
for a non-trivial torus knot is unique.
It seems that such a surface itself is unique up to isotopy in the knot exterior.

\bibliographystyle{amsplain}

%
%
%
\end{document}